\documentclass[10pt]{article}
\usepackage{hyperref}
\usepackage{fancyvrb}  
\usepackage{amsmath,amssymb,color,bbm,ifthen,enumerate}
\usepackage{graphicx}
\usepackage{subfigure}
\usepackage{multirow}
\usepackage{float}

\newtheorem{thm}{Theorem}[section]

\newtheorem{lem}[thm]{Lemma}
\newtheorem{prop}[thm]{Proposition}

\newenvironment{pf}{\paragraph{Proof.}}
{\nopagebreak\hfill\nopagebreak\rule{2mm}{2mm}\par\bigskip}

\begin{document}

\title{ Berry-Esseen and Edgeworth approximations for the tail of an infinite sum of weighted
  gamma random variables}

\author{Mark S. Veillette and Murad S. Taqqu  \thanks{ This work was partially supported
    by the NSF grants DMS-0706786 and DMS-1007616 at Boston University.} 
 }

\maketitle
\begin{abstract}
Consider the sum $Z = \sum_{n=1}^\infty \lambda_n (\eta_n - \mathbb{E}\eta_n)$,
where $\eta_n$ are i.i.d.~gamma random variables with shape parameter $r
> 0$, and the $\lambda_n$'s are
predetermined weights.  We study the asymptotic behavior of the tail
$\sum_{n=M}^\infty \lambda_n (\eta_n - \mathbb{E}\eta_n)$ which is
asymptotically normal under certain conditions. We derive a
Berry-Essen bound and Edgeworth expansions for its
distribution function.  We illustrate the effectiveness of these
expansions on an infinite sum of weighted chi-squared distributions. 

\end{abstract}

\section{Introduction}

Consider a random variable given in terms of an
infinite sum: $\sum_{n=1}^\infty \lambda_n (\eta_n - \mathbb{E}\eta_n)$,
where $\eta_n$ are i.i.d.~gamma random variables with mean
$\mathbb{E}\eta_n = r \theta$ and variance $\mathrm{Var} \  \eta_n = r
\theta^2$, where $r>0$ and $\theta > 0$ are the shape and scale
parameters, respectively.  We may suppose without loss of generality
that $\theta = r^{-1}$, by incorporating the extra parameter into the
constants $\lambda_n$.  We thus consider   
\begin{equation}\label{e:infsumchi2}
Z = \sum_{n=1}^\infty \lambda_n (\eta_n - 1),
\end{equation}
where $\eta_n$ are i.i.d.~gamma with pdf
\begin{equation}\label{e:gamPDF}
f_{\eta_n}(x) = \frac{r^r}{\Gamma(r)} x^{r-1} e^{-r x}, \quad x > 0
\end{equation}
where $r>0$.  We suppose that $\{\lambda_n\}$
is a non-increasing sequence of positive numbers such that $\sum
\lambda_n^2 < \infty$ and are normalized so that
\begin{equation}\label{e:lambdasnormed}
\frac{1}{r} \sum \lambda_n^2 = 1.
\end{equation}
With this setup,
$Z$ has mean zero and variance
\begin{equation}
\mathrm{Var} Z = \sum_{n=1}^\infty \lambda_n^2 \mathrm{Var}
(\eta_n- 1) =  \sum_{n=1}^\infty  \frac{\lambda_n^2}{r} = 1.
\end{equation}
Of particular interest is the case when $\lambda_n =
n^{-\gamma} \ell(n)$, where $\gamma > 1/2$ and $\ell$ is slowly
varying as $n \rightarrow \infty$.  The restriction $\gamma>1/2$
ensures $\sum \lambda_n^2 < \infty$ but allows for cases when either
$\sum \lambda_n = \infty$ or $\sum \lambda_n < \infty$.  

Random variables of the from (\ref{e:infsumchi2}) make up a rich class
of distributions.  Indeed, consider the double Weiner-It\^{o} integral 
\begin{equation}\label{e:WIintegral}
I = \int_{\mathbb{R}^2} '' H(x,y) Z(dx) Z(dy)
\end{equation}
where $Z$ is a complex-valued Gaussian random measure.  
The double prime on the
integral indicates that one excludes the diagonals $\{x = \pm y\}$ from the integration (for more on
integrals of this type, see \cite{taqqu:2010}). In \cite{dobrushin:1979}, Proposition
2, Dobrushin and Major show that the random variable $I$ can be
expressed in the form (\ref{e:infsumchi2}) with $r=1/2$ (chi-squared
distributions).  An important example in this case is the {\it Rosenblatt
  distribution},  discovered by M. Rosenblatt in
\cite{rosenblatt:1961}, and later named after him in
\cite{taqqu:1975}.   For an overview, see \cite{taqqu:2011}.
Properties of the Rosenblatt distribution are further developed in
\cite{veillette:2010a}  using the results we obtain in the present
paper. 

A major difficulty that arises with distributions like (\ref{e:infsumchi2})
is that there is no closed form for its distribution function or
density function.  To make matters worse, even the characteristic
function of $Z$ is not easy to express or compute numerically.   An initial approach to this problem might be to  truncate the sum
(\ref{e:infsumchi2}) at a level $M \geq 1$, and write
$Z = X_M + Y_M$  where
\begin{equation}\label{e:breakupZ}
X_M = \sum_{n=1}^{M-1} \lambda_n(\eta_n - 1), \qquad Y_M =
\sum_{n=M}^\infty \lambda_n(\eta_n - 1)
\end{equation}
and using $X_M$ as an approximation of $Z$ since it is a finite sum of
weighted gamma distributions (an efficient method for computing the PDF/CDF of
such a distribution can be found in \cite{veillette:2010}).  How good is this
approximation?  This question can be partially answered by looking at
the variance of $Y_M$,
\begin{equation}\label{e:defofsigM}
\sigma_M^2 \equiv \mathrm{Var} Y_M =\frac{1}{r}  \sum_{n=M}^\infty \lambda_n^2.
\end{equation}
Depending on the decay of $\lambda_n$, this can tend to $0$ slowly.
For instance, if $\lambda_n \sim C n^{-\gamma}$ for some $\gamma >
1/2$, then
\begin{equation*}
\sigma_M^2 \sim \frac{C}{r}  \int_M^\infty x^{-2\gamma} dx \sim 
\frac{C}{r} M^{1 - 2 \gamma },
\end{equation*}
which tends to $0$ slowly when $\gamma$ is
close to $1/2$, and thus in these cases
$M$ would have to be taken very large for $X_M$ to be a reasonable
approximation.  

Instead of approximating $Z$ by only $X_M$ for $M$ large, we will
instead show when $Y_M$ is asymptotically normal using a Berry-Essen
estimate, and then we will give an
Edgeworth expansion for the distribution function of $Y_M$.  Combining
this with the distribution of $X_M$ will provide a method for
computing the distribution function of $Z$.  This fact can also be
used for simulation of the random variable $Z$ by simulating $X_M$
exactly, and approximating the error with a $\mathrm{N}(0,\sigma_M^2)$
random variable.  

This paper is organized as follows.  In Section \ref{s:LKform}, we
give the characteristic function of $Z$ and $Y_M$ in
L\'{e}vy-Khintchine form.   We then use this form of the
characteristic function to show $Y_M$ is asymptotically normal in
Section \ref{s:BerryEsseen}.  To approximate the CDF of $Y_M$, we prove an approximation lemma in Section
\ref{s:ApproximationLemmaEW} and in Section \ref{s:edgeworth}, we give
an Edgeworth expansion.  Finally, we
demonstrate the accuracy of these approximations in Section
\ref{s:numericalEW} on an example where the $\eta_n$ are chi-squared,
and the sequence $\lambda_n$ is given.

\section{L\'{e}vy-Khintchine representation}\label{s:LKform}

Recall that a random variable $X$ is infinitely divisible if for any
positive integer $n$, one can find i.i.d.~random variables
$X_{1,n},X_{2,n},\dots,X_{n,n}$ such that 
\begin{equation*}
X  \overset{d}{=} X_{1,n} + X_{2,n} + \dots + X_{n,n}
\end{equation*}
The characteristic function of any real valued infinitely divisible random
variable $X$ with $\mathbb{E} X^2 < \infty$ can be expressed in the
following form, known as the L\'{e}vy-Khintchine form. 
\begin{equation}\label{e:LKformulaBS}
\mathbb{E}e^{i u X} = \exp\left( i a u - \frac{1}{2} u^2 \sigma^2 + 
  \int_{\mathbb{R} \setminus \{0\} } (e^{i u x} - 1 -  i u x ) \Pi(dx) \right)
\end{equation}
where $a \in \mathbb{R}$, $\sigma^2 > 0$ and $\Pi$ is a measure on $\mathbb{R}
\setminus \{0\}$, known as the {\it L\'{e}vy measure}, which
satisfies 
\begin{equation}\label{e:LevyCond}
\int_{\mathbb{R} \setminus \{0\}} \min(x^2,1) \Pi(dx) < \infty.
\end{equation}  
For background on such distributions see
\cite{steutel:2004},\cite{sato:1999}, \cite{Bertoin:1996}, or
\cite{Applebaum:2004}.

The random variable $\eta $ with PDF (\ref{e:gamPDF})
is infinitely divisible and has characteristic function
\begin{equation}
\mathbb{E} e^{i u \eta} = \exp\left( \int_0^\infty (e^{i u x} - 1 )
  \Pi(x) \right)
\end{equation}
where the
L\'{e}vy measure is given by $\Pi(dx) = r x^{-1} e^{- r x}  dx$ for $x > 0$
(\cite{Applebaum:2004}, example 1.3.22).       
Hence, if $\lambda > 0$,  the random variable $\lambda(\eta - 1)$ is also infinitely divisible
and its characteristic function is given by
\begin{equation}
\mathbb{E} \exp\left( i u \lambda (\eta- 1) \right) = \exp\left(
  \int_0^\infty (e^{i u x} - 1 - i u x) \left(
    \frac{r}{x} \exp\left( -\frac{r x}{\lambda} \right)   \right) dx \right).
\end{equation}
 By taking an infinite sum of
such distributions as in (\ref{e:infsumchi2}), it is not surprising
that the resulting distribution is also infinitely divisible as
indicated in the following proposition.  

\begin{prop}\label{p:LKFormsum}
The characteristic function of $Z$ defined in (\ref{e:infsumchi2}) is
given by
\begin{equation}\label{e:LKformROS1}
\mathbb{E}e^{i u Z} = \exp\left( \int_0^\infty (e^{i u x} - 1 - i u x)
  \nu(x) dx \right) 
\end{equation}
where $\nu$ is defined as
\begin{equation}\label{e:defofnu}
\nu(x) \equiv \frac{r}{x} \sum_{n=1}^\infty \exp\left( -\frac{r x}{
    \lambda_n} \right).
\end{equation}
\end{prop} 
\begin{pf}
We have
\begin{eqnarray}
\mathbb{E}e^{i u Z} &=& \lim_{M \rightarrow \infty} \mathbb{E} e^{i u
  X_M} \nonumber\\
&=& \mathbb{E} \exp\left( \lim_{M \rightarrow \infty}\int_0^\infty (e^{i u x} -
  1 - i u x) \left( \frac{r}{x } \sum_{n=1}^{M-1}  \exp\left(
      -\frac{r x}{ \lambda_n}
    \right)  \right)dx \right) .
\end{eqnarray}
To pass the limit through the integral above, note that  $|e^{i u x} -
1 - i u x| \leq \frac{1}{2} u^2 x^2$ and thus it suffices to show (using the
dominated convergence theorem) that
\begin{equation}\label{e:intofLMEW}
\frac{r}{2} \int_0^\infty x \sum_{n=1}^\infty \exp\left( -\frac{r x}{
    \lambda_n} \right) dx < \infty.
\end{equation}
This follows since 
\begin{eqnarray*}
\int_0^\infty x \sum_{n=1}^\infty \exp\left( -\frac{r x}{
    \lambda_n} \right) dx = 
\sum_{n=1}^\infty \int_0^\infty x \exp\left( -\frac{r x}{
    \lambda_n} \right) dx = \sum_{n=1}^\infty  \frac{\lambda_n^2}{r^2} < \infty.
\end{eqnarray*}
Thus (\ref{e:intofLMEW}) holds and hence the L\'{e}vy measure is
given by (\ref{e:defofnu}).

\end{pf}

The form (\ref{e:LKformROS1}) of the characteristic function will be useful when we
study the Edgeworth expansion of the tail $Y_M$ defined in
(\ref{e:breakupZ}), whose L\'{e}vy measure is given by
$\nu^{(M)}(x) = r x^{-1} \sum_{n=M}^\infty  \exp\left( -
  \frac{r x}{ \lambda_n} \right)$.

\section{Berry-Esseen Bound}\label{s:BerryEsseen}

In this section we show that under certain conditions on the sequence
$\lambda_n$, then the distribution of the tail $Y_M$ is asymptotically
normal as $M \rightarrow \infty$.  A Berry-Esseen type bound for
infinitely divisible random variables was studied in
\cite{asmussen:2001}, and we will apply a similar method to the random
variable $Y_M$.

Consider the normalized distribution
$\widetilde{Y}_M = \sigma^{-1}_M Y_M$
where $\sigma_M$ is defined in (\ref{e:defofsigM})
and let 
\begin{equation}\label{e:newLM}
\widetilde{\nu}^{(M)}(x) = \sigma_{M}
\nu^{(M)}(\sigma_M x) = \frac{r}{x } \sum_{n=M}^\infty  \exp\left(-
  \frac{r  x \sigma_M}{ \lambda_n} \right), \quad x > 0
\end{equation}
be the density of the L\'{e}vy measure of $\widetilde{Y}_M$.   As the
remark below indicates, $\widetilde{Y}_M = \sigma_M^{-1} Y_M$ does not
always converge to a normal distribution.  
To determine whether $\tilde{Y}_M$ is asymptotically normal, it suffices to
consider the third cumulant of $\tilde{Y}_M$ which we denote by
\begin{equation}\label{e:defofkap3}
\kappa_{3,M} \equiv \int_0^\infty x^3 \widetilde{\nu}^{(M)}(x) dx =
r \sum_{n=M}^\infty \int_0^\infty x^2 e^{-r x \sigma_M/
  \lambda_n} dx = 2 r^{-2} \sigma_M^{-3} \sum_{n=M}^\infty \lambda_n^3.
\end{equation}
The following theorem uses a Berry-Esseen bound to show $\tilde{Y}_M$ is asymptotically normal
if $\kappa_{3,M} \rightarrow 0$.  The constant $0.7056$ appearing in
this bound is the smallest known to date, see \cite{shevtsova:2006}.

\begin{thm}
Let $Z$ be given by (\ref{e:infsumchi2}) and suppose the sequence
$\lambda_n$ is such that
\begin{equation}\label{e:BerEssCond}
\frac{\displaystyle \sum_{n=M}^\infty \lambda_n^3}{\displaystyle
  \left(\sum_{n=M}^\infty \lambda_n^2 \right)^{3/2}} \longrightarrow 0
\quad \mbox{as} \ M\rightarrow \infty.
\end{equation}
Then, $\tilde{Y}_M \rightarrow \mathrm{N}(0,1)$ as $M\rightarrow \infty$
and we have
\begin{equation}\label{e:berryessen}
\sup_{x \in \mathbb{R}}  \left| \mathrm{P}[\tilde{Y}_M \leq x] -
  \Phi(x) \right| \leq 0.7056 \
\kappa_{3,M}
\end{equation}
where $\Phi$ is the standard normal CDF and $\kappa_{3,M}$ is defined in (\ref{e:defofkap3}).
\end{thm}

\noindent {\bf Remark:} It can easily be checked that condition (\ref{e:BerEssCond}) is satisfied
if $\lambda_n$ decays as a power law, i.e.~if $\lambda_n \sim C n^{-\gamma}$ for some $\gamma >
1/2$.  However (\ref{e:BerEssCond}) is {\it not} satisfied if $\lambda_n$
decays exponentially, and in this case convergence to
$\mathrm{N}(0,1)$ will not always hold.  For example, suppose $\lambda_n = 2^{-n-1}$.
Then $Y_M = \sum \lambda_n (\eta_n - 1) = \sum \lambda_n
\eta_n  - \sum \lambda_n \geq -\sum \lambda_n$ and so 
\begin{equation*}
\sigma_M^{-1} Y_M \geq \left( r^{-1} \sum_{n=M}^\infty \lambda_n^2
\right)^{-1/2} \left( -\sum_{n=M}^\infty \lambda_n\right) = -\left( \sum_{n=M}^\infty 2^{-n-1} \right)
\big \slash  \sqrt{r^{-1} \sum_{n=M}^\infty 2^{-2n-2} } = -\sqrt{3r}
\end{equation*}
for all $M$.  Since the normalized random variable $\sigma_M^{-1}Y_M$
is bounded below, it cannot converge in
distribution to $\mathrm{N}(0,1)$.

\begin{pf}
Since $Y_M$ is infinitely divisible, for each $n \geq 1$ we have
\begin{equation} \label{e:iidsum}
Y_M = Y_{M,1}^{(n)} +  Y_{M,2}^{(n)} + \dots + Y_{M,n}^{(n)},
\end{equation}
where $ Y_{M,i}^{(n)}$, $i=1,2,\dots,n$ are i.i.d.~with
mean 0 and variance $\sigma_M^2/n$.  Applying the Berry-Esseen Theorem  (\cite{Gut:2005}, Theorem 7.6.1) 
to the sum (\ref{e:iidsum}),  we have for any $n \geq 1$, 
\begin{eqnarray}
\sup_{x \in \mathbb{R}} \left|\mathrm{P} \left[ \sigma_M^{-1} Y_M \leq x \right] -
\Phi(x) \right| &=& \sup_{x \in \mathbb{R}} \left| \mathrm{P} \left[\frac{1}{\sigma_M
  \sqrt{n}}  (\sqrt{n} Y_M) \leq x \right] -
\Phi(x) \right|\nonumber \\  
&\leq&  0.7056 \frac{ \mathbb{E} \left[ \sqrt{n} Y_{M,1}^{(n)} \right] ^3}{
  \sigma_M^3 \sqrt{n} } \nonumber\\
&=&  0.7056 \frac{ n \mathbb{E} [Y_{M,1}^{(n)}]^3}{
    \sigma_M^3 } .\label{e:nfixed}
\end{eqnarray}
Using Lemma 3.1 in \cite{asmussen:2001}, $n \mathbb{E}[\sigma_M^{-1}  Y_{M,1}^{(n)}]^3
\rightarrow \int_0^\infty x^3 \tilde{\nu}_D^{(M)} (x) dx = \kappa_{3,M}$.  
Thus, we let $n \rightarrow \infty$ in (\ref{e:nfixed}), which gives
(\ref{e:berryessen}).  

To see that the right hand side of this bound tends to $0$ as $M
\rightarrow \infty$, notice that
by (\ref{e:defofsigM}),
\begin{equation}
\kappa_{3,M} = 2 r^{-2}  \sigma_M^{-3}  \sum_{n=M}^\infty \lambda_n^3
=  2 r^{-2}
\left( \frac{1}{r} \sum_{n=M}^\infty   \lambda_n^2 \right)^{-3/2}  \sum_{n=M}^\infty
\lambda_n^3 ,
\end{equation}
which tends to $0$  by the assumption (\ref{e:BerEssCond}), implying
convergence to $\mathrm{N}(0,1)$.  This finishes the proof.  

\end{pf}

If $\lambda_n \sim C n^{-\gamma}$
for $\gamma > 1/2$, then $\kappa_{3,M} = O(M^{-1/2})$ (see
(\ref{e:asyofkapM}) below), which describes
the rate at which the right hand side of (\ref{e:berryessen}) tends
to $0$.  While it is nice to have a practical bound on the error made
when approximating the CDF of $Y_M$ with that of a normal, this rate of
convergence may be too slow.  In the next section, we improve this
approximation by using Edgeworth expansions.  These will do a better
job of approximating the CDF of $Y_M$ for small $M$, however it will
no longer be easy to bound the error made in this approximation exactly.

\section{An approximation lemma}\label{s:ApproximationLemmaEW}

The previous section showed that the tail $Y_M$ can be approximated by
a normal distribution for large $M$.   We shall improve the
approximation to the CDF of $Y_M$ using an Edgeworth expansion.  To
establish the Edgeworth expansion, we will need a lemma
involving an approximation of the characteristic function of
$\tilde{Y}_M$ by a polynomial involving the cumulants.  

In the following sections, we will make the following assumption about
the sequence $\lambda_n$:
\begin{equation}\label{e:asyofLambda}
\lambda_n = \ell(n) n^{-\gamma}
\end{equation}
where $\gamma > 1/2$, and $\ell$ is a slowly varying function at
$\infty$.  With this assumption, (\ref{e:lambdasnormed}) is satisfied and
 \begin{equation}\label{e:asyofsigM}
\sigma_M^2 \sim \frac{1}{r} \int_M^\infty \ell(n)^2 n^{-2\gamma} dn =
\frac{\ell(M)^2 M^{1 - 2 \gamma}}{r( 1 - 2\gamma)}.
\end{equation}
Extending the definition of $\kappa_{3,M}$ in
(\ref{e:defofkap3}), we will denote all cumulants of $\tilde{Y}_M$ by
(see \cite{steutel:2004}, Theorem 7.4),
\begin{eqnarray}
\kappa_{k,M} = \int_0^\infty x^k \widetilde{\nu}^{(M)}(x) dx &=&
r  \sum_{n=M}^\infty  \int_0^\infty x^{k-1} \exp\left(-
  \frac{r x \sigma_M}{\lambda_n} \right) dx   \nonumber \\
&=& \frac{(k-1)!}{r^{k-1} \sigma_M^{k} } \sum_{n=M}^\infty \lambda_n^k, \quad k
\geq 2. \label{e:kapMdef}
\end{eqnarray}
Observe that $\kappa_{2,M} = 1$ and as $M \rightarrow \infty$, (\ref{e:asyofLambda}),
(\ref{e:asyofsigM}), (\ref{e:kapMdef}) and properties of slowly varying functions imply
\begin{eqnarray}
\kappa_{k,M} \sim \frac{(k-1)!}{r^{k-1}} \sigma_M^{-k} \int_M^\infty
\ell(n)^k n^{-k\gamma} dn &\sim& C_k (\ell(M)^{-k} M^{-(\frac{k}{2} - k \gamma)})
(\ell(M)^{k} M^{1 - k \gamma} )\nonumber \\
&=& C_k M^{1 - \frac{k}{2}}, \quad k \geq 2  \label{e:asyofkapM}
\end{eqnarray}
for a constant $C_k$.  Notice that in particular, if $k = 3$, then
$\kappa_{3,M} \sim C_3 M^{-1/2}$, which implies
condition (\ref{e:BerEssCond}).  

In view of Proposition \ref{p:LKFormsum}, the difference between the
log of the characteristic function
of $\tilde{Y}_M$ and that of a standard normal is given by the
following function $I_M$ defined as
\begin{equation}\label{e:IM}
I_M(u) \equiv \int_0^\infty \left(e^{i u x} - 1 - i u x -\frac{ (i ux
  )^2}{2} \right) \tilde{\nu}^{(M)}(x) dx,
\end{equation}
which can be rewritten as
\begin{equation}
I_M(u) = \int_0^\infty (e^{i u x} - 1 - i ux) \tilde{\nu}^{(M)}(x) dx -
\left( -\frac{u^2}{2} \right), \label{e:rewriteIM}
\end{equation}
since $\int_0^\infty x^2 \tilde{\nu}^{(M)}(x) dx = \kappa_{2,M} =  1$.
A key step in developing an Edgeworth expansion is approximating the
function $e^{I_M(u)}$ by a polynomial involving the cumulants, which
is done in the following lemma. 
\begin{lem}\label{l:IMbound}
For $N \geq 3$ and $u>0$, we have as $M \rightarrow \infty$,
\begin{equation}\label{e:FDBexpan}
\left| e^{I_M(u)}  - \left[1 + \sum_{ \eta(N)  } \left( \prod_{m=3}^N  \frac{1}{
  k_m!} \left( \frac{ (i u)^m }{m!} \kappa_{m,M} \right)^{k_m}
\right) \right] \right| \leq 
Q_{N}(u) + \frac{u^{3N-3}}{(3!)^{N-1} (N-1)!} \kappa_{3,M}^{N-1}\exp\left(
  \frac{u^3}{6} \kappa_{3,M} \right),
\end{equation}
where  $\eta(N)$ denotes all non-negative indices
$k_3,k_4,\dots,k_N$ such that 
\begin{equation}\label{e:etainLem}
1 \leq k_3 + 2 k_4 + \dots (N-2) k_N \leq N-2
\end{equation} and  $|Q_{N}(u)|$ is bounded by a polynomial in $u$ whose coefficients are
$O\left(M^{-\left( \frac{N-1}{2} \right) } \right)$ as  $M\rightarrow \infty$.
\end{lem}
\noindent {\bf Remark:}  This bound is a complicated function of $u$, but this will cause no
problem because in the proof of Theorem \ref{t:EWexpansion} below, this
bound is multiplied by $e^{-u^2/2}$ and integrated over $u \in [0,\kappa_{3,M}^{-1}]$.

\begin{pf}
By using Taylor's Theorem on the function $e^{i u x}$ for $u \geq 0$, we have for
each $N\geq 2$,
\begin{equation}
I_M(u) = \int_0^\infty \left(e^{i ux} - 1 - iux - \frac{(i u x)^2}{2}
\right) \tilde{\nu}^{(M)}(dx) =\int_0^\infty \left( \sum_{m=3}^N \frac{ (i u
    x)^m }{m!}  + R_N(u x)   \right) \tilde{\nu}^{(M)}(dx)
\end{equation}
where $R_N$ is a remainder which satisfies 
\begin{equation*}
|R_N(u x)| \leq \frac{( u x)^{N+1}}{(N+1)!}. 
\end{equation*}
Using the definition (\ref{e:kapMdef}) of $\kappa_{k,M}$, $I_M$ becomes
\begin{equation}\label{e:Iexpand}
I_M(u) = \sum_{m=3}^N \frac{ (i u)^m }{m!} \kappa_{m,M} +
\widetilde{R}_N(u)
\end{equation}
where now,
\begin{equation}\label{e:Ibound}
|\widetilde{R}_N(u)| \leq \int_0^\infty \frac{(ux)^{N+1}}{(N+1)!}
\tilde{\nu}^{(M)}(x) dx =   \frac{u^{N+1}}{(N+1)!} \kappa_{N+1,M}.
\end{equation}
Notice that $I_M(u) = \widetilde{R}_2(u)$, which follows from
(\ref{e:Ibound}) by setting $N=2$.  

Turning now to $\exp(I_M(u))$, we apply the classical inequality
\begin{equation}\label{e:classicalETOX}
\left| e^z - \sum_{n=0}^r \frac{z^n}{n!}\right| \leq
\frac{z^{r+1}}{r!} e^{|z|}, \quad z \in \mathbb{R}, \ r \geq 0
\end{equation}
to $\exp(I_M(u))$ and using (\ref{e:Ibound}), we get
\begin{eqnarray}
\left| \exp(I_M(u)) - \sum_{n=0}^{N-2} \frac{I_M(u)^n}{n!}  \right|
 &\leq & \frac{|I_M(u)|^{N-1}}{(N-1)!} \exp(|I_M(u)|) \nonumber\\ &=&
\frac{|\widetilde{R}_2(u)|^{N-1}}{(N-1)!} \exp(|\widetilde{R}_2(u)|) \nonumber\\
&\leq& \frac{ u^{3(N-1)} }{ (3!)^{N-1} (N-1)! } \kappa_{3,M}^{N-1}
\exp\left( \frac{u^3}{3!} \kappa_{3,M} \right). \label{e:firstpIbound}
\end{eqnarray}  
Thus, by adding
and subtracting $\sum_{n=1}^{N-2} \frac{I_M(u)^n}{n!}  $ on
the left hand side of (\ref{e:FDBexpan}), we have
\begin{align}
&\left| \exp\left( I_M(u) \right) - \left(1 + \left[ \sum_{ \eta(N)  } \left( \prod_{m=3}^N  \frac{1}{
  k_m!} \left( \frac{ (i u)^m }{m!} \kappa_{m,M} \right)^{k_m}
\right) \right] \right) \right|  \nonumber\\
&\leq \left| \exp\left( I_M(u) \right)  - \sum_{n=0}^{N-2}
  \frac{I_M(u)^n}{n!}  \right| \nonumber \\
&\quad + \left| \sum_{n=0}^{N-2} \frac{I_M(u)^n}{n!}  -  \left( \left[ \sum_{ \eta(N)  } \left( \prod_{m=3}^N  \frac{1}{
  k_m!} \left( \frac{ (i u)^m }{m!} \kappa_{m,M} \right)^{k_m}
\right) \right] \right)  \right|. \label{e:addsub}
\end{align}
Notice (\ref{e:firstpIbound}) gives a bound for the first term in
(\ref{e:addsub}). Thus, to finish the proof it remains to bound the second term in (\ref{e:addsub}).  
To do this, fix $1 \leq n \leq N-2$ and observe that (\ref{e:Iexpand}) implies 
\begin{equation}
\frac{I_M(u)^n}{n!} = \frac{1}{n!} \left( \sum_{m=3}^{N} \frac{ (i u)^m}{m!}
  \kappa_{m,M} + \widetilde{R}_N(u) \right)^n
\end{equation}
Applying the multinomial theorem,  this becomes
\begin{eqnarray}
\frac{I_M(u)^n}{n!} &=& \frac{1}{n!} \sum_{ \{k_m \}_n  } { n \choose k_3,
  k_4, \dots k_N,  k_{N+1} }   \left[ \prod_{m=3}^N 
\left( \frac{ (i u)^m }{m!} \kappa_{m,M} \right)^{k_m}  \right]
\widetilde{R}_N^{k_{N+1}} \nonumber \\
&=&  \sum_{ \{k_m \}_n } \left[ \prod_{m=3}^N \frac{1}{k_m!}
\left( \frac{ (i u)^m }{m!} \kappa_{m,M} \right)^{k_m}  \right]
\frac{\widetilde{R}_N^{k_{N+1}} }{k_{N+1}!}\label{e:aftermultThm}
\end{eqnarray}
where $\{k_m\}_n$ denotes all sets of non-negative integers $k_m$, $3
\leq m \leq N+1$ such that $k_3 + k_4 + \dots + k_N + k_{N+1} =
n$.  By (\ref{e:asyofkapM}), $\kappa_{m,M} = O(1)$.  Moreover, by (\ref{e:Ibound}) and (\ref{e:asyofkapM}), $|\widetilde{R}_{N}|
\leq \frac{u^{N+1}}{(N+1)!} \kappa_{N+1,M} \sim
\frac{u^{N+1}}{(N+1)!}  C_{N+1}
M^{-(N-1)/2}$, thus any term in (\ref{e:aftermultThm}) involving
$\tilde{R}_N$ (that is with $k_{N+1} \geq 1$) can be
grouped into a function $Q^{(1)}_{n,N}(u)$ which is bounded by a
polynomial with positive coefficients which are
$O(M^{-(N-1)/2})$.  Doing this, (\ref{e:aftermultThm}) becomes
\begin{equation}\label{e:aftermultThm1}
\frac{I_M(u)^n}{n!} =  \sum_{ \{k_m \}' _n  } \left[ \prod_{m=3}^N \frac{1}{k_m!}
\left( \frac{ (i u)^m }{m!} \kappa_{m,M} \right)^{k_m}  \right] + Q_{n,N}^{(1)}(u)
\end{equation}
where $\{k_m \}' _n $ denotes all $k_m$, $3 \leq m \leq N$ such that $k_3
+ k_4 + \dots k_N = n$. In the remaining sum, the coefficients are
\begin{equation}
\prod_{m=3}^N \frac{1}{k_m! (m!)^{k_m}} \kappa_{m,M}^{k_m}.
\end{equation}  
Using (\ref{e:asyofkapM}) again, these coefficients are of the order
\begin{eqnarray}
\prod_{m=1}^N \kappa_{m,M}^{k_m} &=& O\left( M^{ -\sum_{m=3}^N  k_m
   \left( \frac{m - 2}{2} \right)} \right) \nonumber\\
&=& O\left( M^{ -\frac{1}{2} \left[
   \sum_{m=3}^N  m k_m - 2 \sum_{m=3}^N k_m \right]
   } \right) \nonumber \\
&=& O\left( M^{ -\frac{1}{2} \left[
   \sum_{m=3}^N  m k_m - 2n \right]
   } \right).  \label{e:orderofCoefs}
\end{eqnarray}
We shall now isolate the terms in the sum (\ref{e:aftermultThm1}) for which 
\begin{equation}\label{e:kbigsums}
\sum_{m=3}^N m k_m \geq N + 2n - 1.
\end{equation}
They form a polynomial $Q^{(2)}_{n,N}(u)$ whose coefficients by
(\ref{e:orderofCoefs}) are of the order 
\begin{equation*}
O\left( M^{ -\frac{1}{2} \left[
   \sum_{m=3}^N  m k_m - 2n \right]
   } \right)   = O\left(M^{-\frac{1}{2} [N + 2n - 1 - 2n]}\right) =
O\left(M^{-(N-1)/2}\right),
\end{equation*}
where we have used the fact that the $k_m$'s  are chosen to satisfy
(\ref{e:kbigsums}).  Thus,
\begin{equation}
\frac{I_M(u)^n}{n!} =  \sum_{ \{k_m \}'' _n  } \left[ \prod_{m=3}^N \frac{1}{k_m!}
\left( \frac{ (i u)^m }{m!} \kappa_{m,M} \right)^{k_m}  \right] +
Q_{n,N}^{(1)}(u) + Q^{(2)}_{n,N}(u) \label{e:breakupIM}
\end{equation}
where $\{k_m\}'' _n $ denotes all $k_m$, $3 \leq m \leq N$ such that $\sum
k_m = n$ and $\sum m k_m \leq N + 2n - 2$.   Notice that by combining
these two inequalities, the $k_m$'s in this sum also satisfy
\begin{equation}\label{e:simplerforminq}
\sum_{m=3}^N (m-2) k_m \leq (N + 2n - 2) - 2n =  N-2.
\end{equation}

Now, returning to the second term in (\ref{e:addsub}), in light of
(\ref{e:breakupIM}), we have 
\begin{align}
\sum_{n=1}^{N-2} \frac{I_M(u)^n}{n!}  &=  \left( \sum_{n=1} ^{N-2} \left[ \sum_{ \{k_m\}''_n  } \left( \prod_{m=3}^n  \frac{1}{
  k_m!} \left( \frac{ (i u)^m }{m!} \kappa_{m,M} \right)^{k_m}
\right) \right] \right)  + \sum_{n=1}^{N-2} [Q_{n,N}^{(1)}(u) +
Q_{n,N}^{(2)}(u) ]  \nonumber \\
&= \left( \sum_{n=1} ^{N-2} \left[ \sum_{ \{k_m\}''_n  } \left( \prod_{m=3}^n  \frac{1}{
  k_m!} \left( \frac{ (i u)^m }{m!} \kappa_{m,M} \right)^{k_m}
\right) \right] \right)  + Q_N(u) \label{e:rightbefore}
\end{align}
where $Q_N(u) := \sum_{n=1}^{N-2} [Q_{n,N}^{(1)}(u) +
Q_{n,N}^{(2)}(u) ] $ is bounded by a polynomial in $u$ whose
coefficients are $O(M^{-\left(\frac{N-1}{2} \right) })$.  As for the double
sum on the right hand side of (\ref{e:rightbefore}), observe that from
(\ref{e:simplerforminq}), this can be 
rewritten as the (single) sum over all $k_i$, $3 \leq i \leq N$ such that 
\begin{align*}
&1 \leq k_3 + k_4 + \dots k_N \leq N-2 \quad \mbox{and} \\
&k_3 + 2 k_4 + \dots + (N-2) k_N \leq N-2.
\end{align*} 
Since $k_3 + k_4 + \dots +k_N \leq k_3 + 2 k_4 + \dots (N-2) k_N$, these two conditions are satisfied if and only if
\begin{equation*}
1 \leq k_3 + 2 k_4 + \dots + (N-2) k_N \leq N-2,
\end{equation*}
which is the definition of $\eta(N)$ in (\ref{e:etainLem}). Thus, 
\begin{equation}\label{e:piece2IM}
\left|\sum_{n=1}^{N-2} \frac{I_M(u)^n}{n!}  - \sum_{\eta(N)} \prod_{m=3}^N
\frac{1}{k_m!} \left( \frac{(i u)^m}{m!} \kappa_{m,M} \right)^{k_m} \right|
\leq |Q_N(u)|.
\end{equation}
This bounds the second term in (\ref{e:addsub}) and completes the proof.

\end{pf}

\section{Edgeworth expansions}\label{s:edgeworth}

We shall improve the approximation to the CDF of $Y_M$ using and
Edgeworth expansion. A two-term Edgeworth expansion of a general sequence of infinitely
divisible distributions are studied in \cite{lorz:1991}.  We apply a
similar method to our case, but with an Edgewood expansion to any order.    

Given a CDF $F$ of a random variable $X$ and a function $G$ (not necessarily a
CDF) we let $\rho$ denote the supremum norm of the difference $F - G$:
\begin{equation*}
\rho(F,G) = \sup_{x \in \mathbb{R}} |F(x) - G(x)|
\end{equation*}
We can bound $\rho(F,G)$ using the characteristic function of $X$ and
the Fourier-Stieltjes transform of $G$.  This is done in the following
lemma which is proved in \cite{ranga:1976}, Lemma 12.2.  

\begin{lem}\label{l:charbound}
Let $\phi$ be a  characteristic function of a random variable $X$ with CDF
$F$.  Let $G$ be a function for which 
\begin{equation*}
\lim_{x \rightarrow -\infty} G(x) = 0, \quad \lim_{x \rightarrow
  \infty} G(x) = 1, \ \ \mbox{and} \ \  \sup_{x \in \mathbb{R}} |G'(x)| < C,
\end{equation*}
for some constant $C$ and let $g(u) = \int_{\mathbb{R}} e^{i u x}
dG(x)$ be the  Fourier-Stieltjes transform  of $G$.  Furthermore, suppose that
\begin{equation*}
\int_{\mathbb{R}} |x| dF(x) < \infty \ \ \mbox{and} \ \
\int_{\mathbb{R}} |x| dG(x) < \infty.
\end{equation*}
Then for every $U>0$ and $t > t_0$,
\begin{equation}\label{e:boundme}
\rho(F,G) \leq  \frac{1}{4 h(t) - \pi} \int_0^U  |\phi(u) - g(u)|
\frac{du}{u} + 4 t h(t) \frac{C}{U}
\end{equation}
where $h$ and $t_0$ are defined as
\begin{equation*}
h(t) = \int_0^t \frac{\sin^2(x)}{x^2} dx, \ t > 0, \quad \mbox{and}  \quad h(t_0) = \frac{\pi}{4}.
\end{equation*}
\end{lem}
This lemma involves two parameters $t$ and $U$, which must balance each
other (making $U$ large decreases the second term on the right hand
side of (\ref{e:boundme}) and increases the
first, and $t$ has the opposite effect).  In our application, $U$ will
tend to infinity and $t$ will be an unspecified constant.   This lemma
will be used to study the convergence of an Edgewood expansion for $\tilde{Y}_M$.  

We can now state a theorem detailing the convergence rate of an
Edgeworth expansion for the CDF of $\widetilde{Y}_M$ as $M \rightarrow
\infty$.  Recall the {\it
 Hermite polynomials} which can be defined as $H_0(x) = 1$ and
\begin{equation*}
H_{k}(x) = (-1)^k e^{x^2/2} \frac{d^k}{dx^k}  e ^{-x^2/2}, \quad  k \geq 1
,
\end{equation*}
see \cite{kuo:2006}, page 157. A simple induction shows that $H_k$ also satisfies the
recursion formula 
\begin{equation}
H_{k+1}(x) = -e^{x^2/2} \frac{d}{dx}  \left( H_k(x) e^{-x^2/2} \right),  
\quad k \geq 0. \label{e:hermitepoly}
\end{equation}
The first few $H_k$ are given by $H_1(x) = x$, $H_2(x) = x^2 - 1$,
$H_3(x) = x^3-3x, H_4(x) = x^4 - 6x^2 + 3$, $H_5(x) = x^5 - 10x^3 +
15x$, $\dots$.  

The following theorem provides an Edgeworth expansion
of $\tilde{Y}_M$ up to an arbitrary order $N \geq 2$.   

\begin{thm}\label{t:EWexpansion}
As $M \rightarrow \infty$, for each $N \geq 2$ the CDF of $\tilde{Y}_M$ satisfies
\begin{equation}\label{e:EWexpand}
\mathrm{P}[  \tilde{Y}_M \leq x] = \Phi(x) - \phi(x)\left\{
  \sum_{\eta(N)} \left[ \prod_{m=1}^N \frac{1}{k_m!} \left(
      \frac{\kappa_{m,M}}{m!} \right)^{k_m} \right]
  H_{\zeta(k_3,\dots,k_N)} (x) \right\} + O\left( M^{-\frac{N-1}{2}} \right),
\end{equation}
where $\Phi$ and $\phi$ denote the standard normal CDF and PDF, $\eta(N)$ denotes all non-negative indices
$k_3,k_4,\dots,k_N$ such that 
\begin{equation}\label{e:etainThm}
1 \leq k_3 + 2 k_4 + \dots (N-2) k_N
\leq N-2
\end{equation} and 
\begin{equation}\label{e:zetainThm}
\zeta( k_1,\dots,k_N) = 3k_3 + 4 k_4 + \dots N k_N - 1.
\end{equation}
Moreover, the error $O(M^{-(N-1)/2})$ is uniform for all $x \in
\mathbb{R}$.  
\end{thm}

For example, if $N = 2$, there is no solution to (\ref{e:etainThm}).
If $N = 3$, the only solution to  (\ref{e:etainThm}) is $k_3 = 1$.  If
$N=4$, we have the additional solutions $k_3 = 2$, $k_4 = 0$ and $k_3 =
0$, $k_4 = 1$.  Thus, for small values of $N$, the right hand side of (\ref{e:EWexpand}) becomes

\begin{align*}
& \mathbf{N=2:} \quad \mathrm{P}[  \tilde{Y}_M \leq x] = \Phi(x)  +
O(M^{-1/2}) \\
&  \mathbf{N=3:} \quad \mathrm{P}[  \tilde{Y}_M \leq x] = \Phi(x)  -
\phi(x) \frac{ H_2(x)}{3!} \kappa_{3,M}+ O(M^{-1}) \\
& \mathbf{N=4:} \quad \mathrm{P}[  \tilde{Y}_M \leq x] = \Phi(x)  -
\phi(x) \left[ \frac{H_2(x)}{3!} \kappa_{3,M}+ \left( \frac{ H_3(x)}{4!} \kappa_{4,M}
  + \frac{H_5(x)}{(2!)(3!)^2} \kappa_{3,M}^2\right) \right]+
O(M^{-3/2})\\
&\mathbf{N=5:} \quad \mathrm{P}[  \tilde{Y}_M \leq x] = \Phi(x)  -
\phi(x) \left[ \frac{H_2(x)}{3!} \kappa_{3,M}+ \left( \frac{ H_3(x)}{4!} \kappa_{4,M}
  + \frac{H_5(x)}{(2!)(3!)^2} \kappa_{3,M}^2 \right) \right. \\
& \ \  \qquad \qquad \qquad \qquad \qquad \qquad \qquad \left.+ \left( \frac{H_4(x)}{5!} \kappa_{5,M} +
  \frac{H_6(x)}{(3!)(4!)} \kappa_{4,M} \kappa_{3,M} + \frac{ H_8(x) }{
    (3!)^4 } \kappa_{3,M}^3 \right) \right]+ O(M^{-2})
\end{align*}

A more revealing (but slightly more complicated) statement of Theorem
\ref{t:EWexpansion} is
\begin{equation}
\mathrm{P}[  \tilde{Y}_M \leq x] = \Phi(x) - \phi(x)\left\{
  \sum_{n=3}^{N}\left( \sum_{\eta'(n)} \left[ \prod_{m=3}^n \frac{1}{k_m!} \left(
      \frac{\kappa_{m,M}}{m!} \right)^{k_m} \right]
  H_{\zeta(k_3,\dots,k_n)} (x) \right) \right\} + O\left( M^{-\frac{N-1}{2}} \right),
\end{equation}
where $\eta'(n)$ denotes all $k_3, k_4,\dots,k_N$ such that $k_3 + 2
k_4 + \dots + (n-2)k_N = n-2$.  In this form, it is clearer what
additional terms appear in the expansion as you increase $n$ from 3 to
$N$.


\subsubsection*{Proof of Theorem \ref{t:EWexpansion}}

\begin{pf}
 Define $G(x)$ as 
\begin{equation}
G(x) = \Phi(x) - \phi(x) \left\{
  \sum_{\eta(N)} \left[ \prod_{m=1}^N \frac{1}{k_m!} \left(
      \frac{\kappa_{m,M}}{m!} \right)^{k_m} \right]
  H_{\zeta(k_3,\dots,k_N)} (x) \right\}.
\end{equation}
Then by  (\ref{e:hermitepoly}), we also have
\begin{equation*}
\frac{dG}{dx} = \phi(x) \left( 1 + \sum_{\eta(N)} \left[ \prod_{m=1}^N
    \frac{1}{k_m!} \left( \frac{\kappa_{m,M}}{m!} \right)^{k_m} \right]
  H_{\zeta(k_3,\dots,k_N) + 1} \right).
\end{equation*} 
Using this and the fact
$\int_{\mathbb{R}} H_k(x) \phi(x) e^{i u x} dx = (-1)^k \int_{\mathbb{R}}
\left( \frac{d^k}{dx^k} \phi(x) \right) e^{i u x} dx = (i
u)^k e^{-u^2/2}$, the Fourier-Stieltjes transform of $G$ is given by
\begin{eqnarray}
g(u) = \int_{\mathbb{R}} e^{i u x} dG(x) &=& \exp\left(-\frac{u^2}{2}
\right) \left(   1 + \sum_{ \eta(N)  } \left( \prod_{m=3}^N  \frac{1}{
  k_m!} \left( \frac{ \kappa_{m,M} }{m!} \right)^{k_m} (i
u)^{\zeta(k_3,\dots,k_N) + 1}
\right)  \right) \nonumber\\
&=& \exp\left(-\frac{u^2}{2}
\right) \left(   1 + \sum_{ \eta(N)  } \left( \prod_{m=3}^N  \frac{1}{
  k_m!} \left( \frac{ (i u)^m }{m!} \kappa_{m,M} \right)^{k_m}
\right)  \right) \label{e:defofg_Ros}
\end{eqnarray}
where we have used the definition of $\zeta$ in (\ref{e:zetainThm}).
Let $\varphi^{(M)}(u)$ be the characteristic function of
$\widetilde{Y}_M$:
\begin{equation}
\varphi^{(M)}(u) = \exp\left( \int_0^\infty (e^{i u x} - 1 - i u x)
  \widetilde{\nu}^{(M)}(x) dx \right).
\end{equation}

Since $N \geq 2$, choose $\epsilon > 0$ such that
\begin{equation*} 
\kappa_{3,M}^{-1} <  \epsilon
\kappa_{N+1,M}^{-1} 
\end{equation*} 
for all $M \geq 1$ (this exists by
(\ref{e:asyofkapM})).  To show (\ref{e:EWexpand}) using Lemma
\ref{l:charbound}, it suffices to show that
\begin{equation}\label{e:JM}
J_M := \int_0^{U_M} |\varphi^{(M)}(u) - g(u) | \frac{du}{u} = O\left(M^{-(N-1)/2}\right),
\end{equation}
 where 
\begin{equation}\label{e:defofUM_EW}
U_M := \epsilon \kappa_{N+1,M}^{-1} \sim C_{N+1} M^{-\left(1 -
   \frac{N+1}{2}\right)}= C_{N+1} M^{\frac{N-1}{2}} 
\end{equation}
from (\ref{e:asyofkapM}). Notice that with this choice of $U_M$, the second term on the right
hand side of (\ref{e:boundme}) is already of order $O(U_M^{-1}) = O(M^{-(N-1)/2})$ and
thus we need to only bound $J_M$.   

Using (\ref{e:rewriteIM}), notice that
\begin{eqnarray*}
\varphi^{(M)} (u) = \exp\left( \int_0^\infty (e^{i u x} - 1 - i ux)
  \widetilde{\nu}^{(M)}(x) dx \right) = \exp\left( I_M(u) - \frac{u^2}{2} \right).
\end{eqnarray*}
Using this and the definition of $g$ in (\ref{e:defofg_Ros}), we can break up the
 integral $J_M$ in  (\ref{e:JM}) as
\begin{eqnarray*}
J_M &=& \int_0^{U_M} \exp\left( - \frac{u^2}{2} \right) \left(
  \exp(I_M(u)) -  \left(   1 + \sum_{ \eta(N)  } \left( \prod_{m=3}^N  \frac{1}{
  k_m!} \left( \frac{ (i u)^m }{m!} \kappa_{m,M} \right)^{k_m}
\right)  \right) \right)
\frac{du}{u} \\
&:=& J_{M,1} + J_{M,2} + J_{M,3},
\end{eqnarray*}
where
\begin{align*}
& J_{M,1} = \int_0^{\kappa_{3,M}^{-1}} \exp\left( - \frac{u^2}{2} \right) \left(
  \exp(I_M(u)) -   \left(   1 + \sum_{ \eta(N)  } \left( \prod_{m=3}^N  \frac{1}{
  k_m!} \left( \frac{ (i u)^m }{m!} \kappa_{m,M} \right)^{k_m}
\right)  \right) \right)
\frac{du}{u}  \\
& J_{M,2} = \int_{ \kappa_{3,M}^{-1} }^{U_M} \exp\left( -\frac{u^2}{2} + I_M(u)
\right) \frac{du}{u} \\
& J_{M,3} = -\int_{\kappa_{3,M}^{-1}}^{U_M} \exp\left( -\frac{u^2}{2} \right) \left(   1 + \sum_{ \eta(N)  } \left( \prod_{m=3}^N  \frac{1}{
  k_m!} \left( \frac{ (i u)^m }{m!} \kappa_{m,M} \right)^{k_m}
\right)  \right)  \frac{du}{u} .
\end{align*}

We will now show that $J_{M,i} = O(M^{-(N-1)/2})$, $i=1,2,3$, which
with the help of Lemma \ref{l:charbound}, 
will imply the result.

\medskip 

\noindent {\bf Estimate for $J_{M,1}$:}

\medskip

From Lemma \ref{l:IMbound}, we have that
\begin{equation*}
|J_{M,1}| \leq \int_0^{\kappa_{3,M}^{-1}} \exp\left( -\frac{u^2}{2}
\right) \left| Q_{N}(u) + \frac{u^{3N-3}}{(3!)^{N-1} (N-1)!} \kappa_{3,M}^{N-1}\exp\left(
  \frac{u^3}{6} \kappa_{3,M} \right) \right| \frac{du}{u}
\end{equation*}
where $Q_N$ is bounded by a polynomial in $u$ whose coefficients are
$O(M^{-(N-1)/2}$.  Thus,
\begin{equation}
|J_{M,1}| \leq \int_0^{\kappa_{3,M}^{-1}} \exp\left( -\frac{u^2}{2}
\right) \left| Q_{N}(u) \right| \frac{du}{u}  +
\int_0^{\kappa_{3,M}^{-1}} \exp\left( -\frac{u^2}{3} \right)
\frac{u^{3N - 4}}{(3!)^{N-1} (N-1)!} \kappa_{3,m}^{N -1} du, \label{e:JMbound1}
\end{equation}
since on the interval $0 < u < \kappa_{3,M}^{-1}$, we have $u \kappa_{3,M}
\leq 1$ and
 \begin{equation*}
\exp\left( -\frac{u^2}{2}\right) \exp\left( \frac{u^3}{6}
   \kappa_{3,M}\right) \leq
 \exp\left(-\frac{u^2}{3}\right). \end{equation*}
 The first term in (\ref{e:JMbound1}) is $O(M^{-(N-1)/2})$ since all
coefficients of $Q_N$ are of
this order.  The second term in (\ref{e:JMbound1}) is also of order
$O(M^{-(N-1)/2})$ by (\ref{e:asyofkapM}).  Thus,  $J_{M,1} =
O(M^{-(N-1)/2})$.    

\medskip 

\noindent {\bf Estimate for $J_{M,2}$:}

\medskip

We will in fact show that $J_{2,M} = o(M^{-(N-1)/2})$. First, observe that
\begin{equation} \label{e:B1forJ2}
|J_{2,M}| \leq \int_{\kappa_{3,M}^{-1}}^{U_M} \exp\left( 
    -\frac{u^2}{2}  + \mathrm{Re} [I_M(u) ] \right) \frac{du}{u} .
\end{equation}
Thus, we must show the integrand tends to zero fast enough.  Notice,
using (\ref{e:rewriteIM}), that
\begin{eqnarray*}
-\frac{u^2}{2}  + \mathrm{Re} [I_M(u) ]  = -\int_0^\infty ( 1 - \cos( ux) ) \widetilde{\nu}^{(M)}(x) dx
\end{eqnarray*}
Using (\ref{e:newLM}), we compute this integral:
\begin{eqnarray}
 A_M(u) &:=& \int_0^\infty (1 -  \cos( ux) ) \widetilde{\nu}^{(M)}(x) dx
 =  r \int_0^\infty ( 1 - \cos( ux) 
) \left( \sum_{n=M}^\infty \frac{e^{-x r \sigma_M/ \lambda_n
      }}{x} \right)  dx \nonumber\\
 &=& r \sum_{n=M}^\infty \int_0^\infty \frac{ 1 - \cos(ux)  }{ x}
  e^{-x r \sigma_M/ \lambda_n
      } dx = r \sum_{n=M}^\infty  \int_0^\infty \int_0^u \sin(t x)
e^{-x r \sigma_M/\lambda_n }dt dx \nonumber\\
&=&  r \sum_{n=M}^\infty  \int_0^u \left( \int_0^\infty \sin(t x)
e^{-x r \sigma_M/ \lambda_n } dx \right) dt = r \sum_{n=M}^\infty  \int_0^u \left( \frac{t}{t^2 +
    \frac{r^2 \sigma_M^2}{ \lambda_n^2} } \right) dt \nonumber \\
&=& \frac{r}{2}  \sum_{n=M}^\infty \log\left( 1 + \frac{u^2 \lambda_n^2
  }{r^2 \sigma_M^2} \right), \label{e:Lastpartlongint}
\end{eqnarray}
where we have used the integral identity $\int_0^\infty \sin( t x)
e^{-zx} dx = t/(t^2 + z^2) $ in the fourth line, which can be shown by integration
by parts.

Using the properties of slowly varying functions,  for any $\gamma'$
and $\gamma''$ for which $\gamma' >
\gamma > \gamma''$ one can find constants $\alpha_1,\alpha_2$ such that $
\alpha_2 n^{-\gamma''} \gtrsim \lambda_n \gtrsim
\alpha_1 n^{-\gamma'}$, that is, $\limsup \lambda_n
n^{\gamma''}/\alpha_2 \leq 1$ and  $\liminf \lambda_n n^{\gamma'}/\alpha_1
\geq 1$.  Since (\ref{e:Lastpartlongint}) is increasing is
$u$, on the interval $ \kappa_{3,M}^{-1} < u < U_M$, as $M \rightarrow \infty$,
\begin{eqnarray}
A_M(u) \geq A_M(\kappa_{3,M}^{-1}) &=&  \frac{r}{2}  \sum_{n=M}^\infty \log\left( 1 + \frac{   \lambda_n^2
  }{r^2 \kappa_{3,M}^2 \sigma_M^2} \right) \nonumber \\
&\gtrsim& \frac{r}{2} \int_M^\infty \log\left( 1 + \frac{  \alpha_1^2
    y^{-2\gamma'} 
  }{ r^2 \kappa_{3,M}^2 \sigma_M^2} \right) dy
\\
&=& \frac{r}{2} \int_M^\infty \log\left( 1 + (\beta_M
  y)^{-2\gamma'}  \right) dy \label{e:RbeforeCOV}
\end{eqnarray}
where
\begin{equation*}
 \beta_M = \left( \frac{ \alpha_1^2 }{r^2 \kappa_{3,M}^2 \sigma_M^2}
\right)^{-1/2\gamma' }. 
\end{equation*}
Making the change of variables $w = \beta_M y$, the integral (\ref{e:RbeforeCOV}) becomes
\begin{equation*}
\frac{r}{2 \beta_M}   \int_{M \beta_M  }^\infty \log\left( 1 +
  w^{-2\gamma'} \right) dw
\end{equation*}
Equations (\ref{e:asyofsigM}) and (\ref{e:asyofkapM}) together with
the choice of $\gamma''$ imply
$\kappa_{3,M}^2 \sigma_M^2 \lesssim C  M^{-2\gamma''}$ for a constant $C > 0$, hence
\begin{equation*}
\beta_M \lesssim C' M^{-\gamma''/ \gamma'}  
\end{equation*} 
for another constant $C'$.   Thus, we have shown that as $M \rightarrow \infty$,
\begin{equation*}
A_M(\kappa_{3,M}^{-1}) \gtrsim \frac{ r C'}{ 2} M^{\gamma''/\gamma'}
\int_{C' M^{1 - \gamma''/\gamma'}}^\infty \log( 1 + w^{-2\gamma'})
dw .
\end{equation*}
Since $\log(1 + x) \sim x$ as $x \rightarrow 0$, we have
\begin{eqnarray}
 A_M(\kappa_{3,M}^{-1}) &\gtrsim& \frac{r C'}{2 } M^{\gamma''/\gamma'}
 \int_{C' M^{1 - \gamma''/\gamma'}}^\infty w^{-2\gamma'} dw \nonumber\\
&=& C'' M^{\gamma''/\gamma'}
M^{(1-\gamma''/\gamma')(1-2\gamma')} \nonumber \\
&=& C'' M^{1 - 2(\gamma'-\gamma'')} \label{e:asyofA}
\end{eqnarray}
for another constant $C''$.  Notice (\ref{e:asyofA}) tends to infinity
so long as $\gamma'-\gamma''$ is chosen to be smaller than $1/2$.
Now, returning to $J_{2,M}$, (\ref{e:defofUM_EW}), (\ref{e:B1forJ2}), (\ref{e:asyofA}), and (\ref{e:asyofkapM})
imply
\begin{eqnarray*}
|J_{2,M}| &\leq& \int_{\kappa_{3,M}^{-1}}^{U_M} \exp(-A_M(u))
\frac{du}{u} \\
&\leq& \exp(-A_M( \kappa_{3,M}^{-1} ) ) \int_{\kappa_{3,M}^{-1}}^{U_M}
\frac{du}{u} \\
&=& \exp(-A_M(\kappa_{3,M}^{-1})) \log\left(
  \frac{\epsilon \kappa_{3,M}}{\kappa_{N+1,M}} \right) \\
&\lesssim& C''' \exp\left( -C'' M^{1 - 2(\gamma'-\gamma'') } \right)\log M, \quad \mbox{for some $C'''>0$} \\
&=& o(M^{-(N-1)/2}).
\end{eqnarray*}

\medskip 

\noindent {\bf Estimate for $J_{M,3}$:}

\medskip

For $J_{M,3}$, we have
\begin{equation}
|J_{M,3}| \leq \left| \int_{\kappa_{3,M}^{-1}}^{U_M} \exp\left(
    -\frac{u^2}{2} \right) \left(   1 + \sum_{ \eta(N)  } \left( \prod_{m=3}^N  \frac{1}{
  k_m!} \left( \frac{ (i u)^m }{m!} \kappa_{m,M} \right)^{k_m}
\right)  \right)  \frac{du}{u}  \right|. \nonumber
\end{equation}
By bounding all the coefficients of the polynomial in $u$ by their
maximum value, we have
\begin{eqnarray}
|J_{M,3}|  &\leq& \max_{\eta(N)} \left(1, \prod_{m=3}^N \kappa_{m,M}^{k_m} \right)
\int_{\kappa_{3,M}^{-1}}^\infty \exp\left( -\frac{u^2}{2} \right)
\left( 1 +\sum_{ \eta(N)  } \left( \prod_{m=3}^N  \frac{1}{
  k_m!} \left( \frac{ u^m }{m!} \right)^{k_m} \right)  \right) du \nonumber
\\
&\sim&   \max_{\eta(N)} \left(1, \prod_{m=3}^N \kappa_{m,M}^{k_m}
\right)  \int_{\ell(M)^{-1} M^{1/2}}^\infty  \exp\left(
  -\frac{u^2}{2} \right) p(u) du  \label{e:JMbound3}
\end{eqnarray}
where $p(u)$ is a polynomial in $u$ whose coefficients do not depend
on $M$.   Choosing a constant $C>0$ large
enough such that $p(u) \leq C e^u$ for all $u > \ell(M)^{-1} M^{1/2}$,
we see
\begin{eqnarray}
\int_{\ell(M)^{-1} M^{1/2}}^\infty  \exp\left(
  -\frac{u^2}{2} \right) p(u) du  &\leq& \int_{\ell(M)^{-1} M^{1/2}}^\infty  \exp\left(
  -\frac{u^2}{2} \right) Ce^{u}  du \nonumber \\
&=& C e^{1/2} \int_{\ell(M)^{-1} M^{1/2} -1 }^\infty   e^{-w^2/2} dw \\
 &=&       C \sqrt{ \frac{e \pi}{2} }
\mathrm{Erfc}\left( \frac{ \ell(M)^{-1} M^{1/2} - 1}{\sqrt{2}} \right) \label{e:erfcbound}
\end{eqnarray} 
where $\mathrm{Erfc}(u) = \frac{2}{\sqrt{\pi}} \int_u^\infty
e^{-w^2} dw$.  Using the fact that $\mathrm{Erfc}(u) \sim
e^{-u^2}/(\sqrt{\pi} u)$, (\cite{oldham:2009}, equation 40:9:2)   (\ref{e:JMbound3}) and (\ref{e:erfcbound}) imply
\begin{eqnarray*}
|J_{M,3}| &\leq& \max_{\eta(N)} \left(1, \prod_{m=3}^N \kappa_{m,M}^{k_m}
\right)  C \sqrt{ \frac{e \pi}{2} }
\mathrm{Erfc}\left( \frac{ \ell(M)^{-1} M^{1/2} - 1}{\sqrt{2}} \right)
\\
&\sim& O(1)   \frac{\exp\left( -\frac{1}{2} \ell(M)^{-2} M
\right)}{ \ell(M)^{-2} M
} \\
&=& o(M^{-(N-1)/2})
\end{eqnarray*}

Combining the estimates for $J_{M,i}$, $i=1,2,3$, together with
Lemma \ref{l:charbound} implies the desired result.  

\end{pf}

\section{A Numerical example}\label{s:numericalEW}

In this section, we will demonstrate the utility of the Edgeworth
expansion given in Theorem \ref{t:EWexpansion} for computing the CDF of random variable
of the form (\ref{e:infsumchi2}).  Consider the example where $r =
1/2$, i.e.~$\eta_n$ is chi-squared with 1 degree of freedom, and
the $\lambda_n$'s are given simply by
\begin{equation*}
\lambda_n = C n^{-3/4}
\end{equation*}
where the normalization constant is $C = (2 \zeta(6/4) )^{-1/2} =  (2 \sum_{n=1}^\infty n^{-6/4})^{-1/2}
\approx 0.4375$, where $\zeta$ denotes the Riemann zeta function.  To compute the CDF of $Z = \sum_{n=1}^\infty
\lambda_n(\eta_n - 1) = X_M + Y_M$, where $X_M$ and $Y_M$ are defined
in (\ref{e:breakupZ}),  we will proceed in three steps.
\begin{itemize}
\item[1.]  Choose a truncation level $M \geq 1$.  We will see below
  that $M$ does not need to be too large.  Once an $M$ is chosen, one
  must be able to compute the CDF $F_{X_M}(x)$ of $X_M$, which is a finite sum of
  weighted chi-squared random variables. There are multiple techniques
  for doing this, for instance methods based on Laplace transform
  inversion, \cite{veillette:2010}, \cite{castano:2005},  or Fourier
  transform inversion, \cite{abate:1992}.  

\item[2.]  Choose an $N \geq 3$ and compute the appropriate terms in
  Edgeworth expansion for the CDF $F_{Y_M}(x)$ of $Y_M$.  For this
  example, $\sigma_M$ defined in (\ref{e:defofsigM}), and the
  $\kappa_{k,M}$'s defined in (\ref{e:kapMdef}), can be computed in terms of the Riemann Zeta function:
\begin{eqnarray*}
\sigma_M^2 &=& 2 \sum_{n=M}^\infty \lambda_n^2 =   2C^{2} \left( \zeta(2 \gamma) - \sum_{n=1}^{M-1} n^{-6/4} \right)        \\
\kappa_{k,M} &=& 2^{k-1} (k-1)! \sigma_M^{-2} \sum_{n=M}^\infty
\lambda_n^k = 2^{k-1} (k-1)! \sigma_M^{-k} C^{k} \left( \zeta(k \gamma) -
  \sum_{n=1}^{M-1} n^{-3k/4} \right)
\end{eqnarray*}

\item[3.]  The CDF of the sum $Z = X_M + Y_M$ is given by the convolution
\begin{equation}
F_Z(x) = \int_{-\infty}^\infty F_{X_M}(x-y) dF_{Y_M}(y).
\end{equation}
We compute this integral numerically in MATLAB using standard techniques.   
\end{itemize}

\begin{figure}[H]
\begin{center}
\hspace*{-.4in}
\includegraphics[scale=.6]{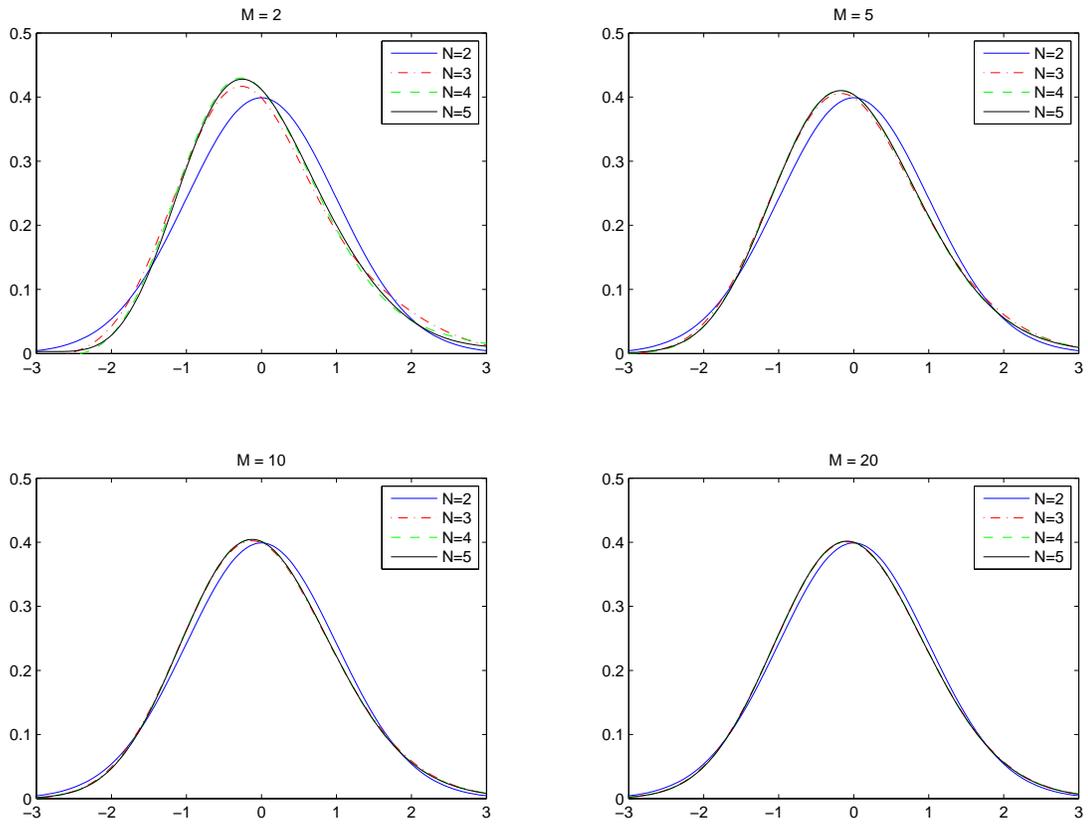}
\caption{ Edgeworth approximations to the density of $\tilde{Y}_M$ for various
values of $M$. We can see that for this example, the converge of $Y_M$
to a normal distribution is fast as $M$ grows and increasing $N$ beyond $5$ causes
a negligible change in the distribution function.} \label{f:EWYM}
\end{center}
\end{figure}

\begin{figure}[H]
\begin{center}
\hspace*{-.8in}
\includegraphics[scale=.4]{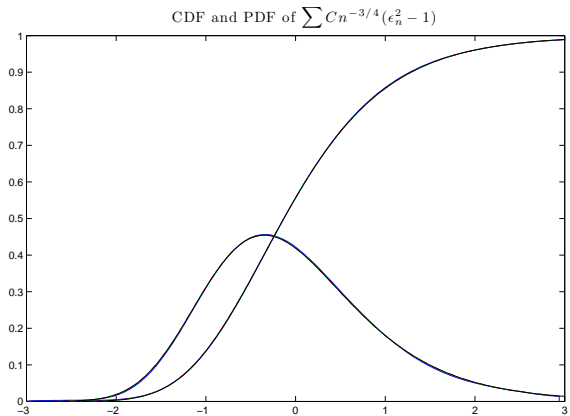}
\caption{ Approximation to the CDF and PDF of $Z$ in the case $\lambda_n = C
  n^{-3/4}$ using the $N=5$ Edgeworth expansion for the tail $Y_M$.
  There are 4 curves corresponding to $M=2,5,10$ and $20$ in both
  curves and are almost indistinguishable suggesting fast convergence of this method.} \label{f:CDFZ}
\end{center}
\end{figure}

We have studied approximations of $Y_M$ for various values of $M$ and
$N$.  Figures \ref{f:EWYM} and \ref{f:CDFZ} give a sense of how good
these approximations are.  We look at the Edgeworth approximations to the density of $Y_M$
and see how these behave as both $M$ and $N$ grow.  
Figure \ref{f:EWYM} shows plots of the $N=2,3,4,5$ Edgeworth approximations to the
density of $\tilde{Y}_M$ for $M=2,5,10$ and $20$. 

An Edgeworth expansion with values of $N \geq 2$ involves corrections
to the normal distribution.  Increasing $N$ improves on this
correction.  If the improvement is already negligible if one goes from
$N=2$ to $N=3$, then the
distribution is close to normal.  This appears to be the case in
Figure \ref{f:EWYM} for already small values of $M$ ($M=10$).  

What happens at smaller values of $M$?  We note that in Figure
\ref{f:EWYM}, that even for $M = 2$, there seems to be no change in
the Edgeworth correction as $N$ goes from $4$ to $5$.  Hence
it appears that for
small values of $M$, a high level of accuracy is already reached by
$N=5$ as it is hard to distinguish the $N=4$ and $N=5$ curves. 

For
this reason, we will use $N=5$ to approximate the CDF of the full
distribution $Z$.  Figure
\ref{f:CDFZ} shows the CDF computed using an $N=5$ Edgeworth expansion
for $Y_M$ for various values of $M$.  Since the resulting
approximation is nearly independent of $M=2,5,10,20$, it is clear that
 the convergence of the Edgeworth expansions is fast for this
 example.  The techniques developed here are used in
 \cite{veillette:2010a} to obtain the numerical evaluation of the CDF
 and PDF of the Rosenblatt distribution. 

\end{document}